\documentclass[reqno]{amsart}

\usepackage{graphicx}
\newtheorem{thm}{Theorem}
\newtheorem{cor}{Corollary}
\newtheorem{lem}[thm]{Lemma}
\newtheorem{prop}[thm]{Proposition}
\theoremstyle{definition}
\newtheorem{defn}[thm]{Definition}
\theoremstyle{remark}
\newtheorem{rem}[thm]{Remark}
\DeclareMathOperator{\RE}{Re}
\DeclareMathOperator{\IM}{Im}
\DeclareMathOperator{\ess}{ess}
\newcommand{\eps}{\varepsilon}
\newcommand{\To}{\longrightarrow}
\newcommand{\h}{\mathcal{H}}
\newcommand{\s}{\mathcal{S}}
\newcommand{\A}{\mathcal{A}}
\newcommand{\J}{\mathcal{J}}
\newcommand{\M}{\mathcal{M}}
\newcommand{\W}{\mathcal{W}}
\newcommand{\X}{\mathcal{X}}
\newcommand{\BOP}{\mathbf{B}}
\newcommand{\BH}{\mathbf{B}(\mathcal{H})}
\newcommand{\KH}{\mathcal{K}(\mathcal{H})}
\newcommand{\Real}{\mathbb{R}}
\newcommand{\Complex}{\mathbb{C}}
\newcommand{\Field}{\mathbb{F}}
\newcommand{\RPlus}{\Real^{+}}
\newcommand{\Polar}{\mathcal{P}_{\s}}
\newcommand{\Poly}{\mathcal{P}(E)}
\newcommand{\EssD}{\mathcal{D}}
\newcommand{\Lom}{\mathcal{L}}
\newcommand{\States}{\mathcal{T}}
\newcommand{\abs}[1]{\left\vert#1\right\vert}
\newcommand{\set}[1]{\left\{#1\right\}}
\newcommand{\seq}[1]{\left<#1\right>}
\newcommand{\norm}[1]{\left\Vert#1\right\Vert}
\newcommand{\essnorm}[1]{\norm{#1}_{\ess}}

\begin{document}

\title[Transport on a rotating plane]
 {The solution of the cauchy problem for the two-dimensional transport equation on a rotating plane}

\author{Olga Rozanova, Olga Uspenskaya}

\address{Department of Mechanics and Mathematics, Lomonosov Moscow State University, Moscow, Russia}

\email{rozanova@mech.math.msu.su}




\dedicatory{}



\begin{abstract}
  The
limiting case of the system of equations of two-dimensional gas
dynamics in the presence of the Coriolis force, which can be
obtained under the assumption of a small pressure, is considered.
With this approach, the equation for the velocity vector (transport
equation) is split off from the system and can be solved separately.
Using the method of stochastic perturbation along characteristics,
we obtain an explicit asymptotic representation of a smooth solution
of transport equations and analyze the process of formation of
singularities of solution using a specific example. It is concluded
that the presence of the Coriolis force prevents the singularities
formation.
\end{abstract}

\maketitle

\section*{Introduction}

The system of equations of dynamics of polytropic gas on a plane in a coordinate system rotating with a constant angular velocity is widely used in meteorology \cite {Vallis}. It has the form:
\begin{equation}\label{e1}
 \rho(\partial_t {\bf U}+ ({\bf U},\nabla) {\bf U} + l L {\bf U})= -\nabla p,
\end{equation}
\begin{equation}\label{e2}
 \partial_t \rho + {\rm div}(\rho {\bf U})=0,
 \end{equation}
 \begin{equation}\label{e3}
 \partial_t S +({\bf U},\nabla S)=0,
\end{equation}
where $\rho\ge 0$, $p\ge 0$, ${\bf U}=(U_1, U_2)$, $S$ are density, pressure, velocity and entropy, respectively.  They are functions of time
$t\ge 0$ and coordinate $x\in {\mathbb R}^2,$ $\quad L = \left(\begin{array}{cr} 0 & -1 \\
1 & 0
\end{array}\right)$, $l=\rm const>0$ is the Coriolis parameter.
The state equation is
\begin{equation}\label{e_state}
p=\rho^\gamma e^S,\qquad \gamma\in (1, 2].
\end{equation}
Equations of atmospheric dynamics of a medium scale are reduced to systems of this kind after some vertical averaging \cite{Obukhov}, which is possible due to the smallness of the vertical scale compared to the horizontal one. Therefore, the behavior of the solutions of such systems is of great interest.
Since the system (\ref{e1}) - (\ref{e_state}) is strictly hyperbolic for positive density, the solution of the Cauchy problem
\begin{equation}\label{CD}
(\rho(t,x), p(t,x), {\bf U}(t,x))|_{t=0} = (\rho^0(x), p^0(x), {\bf
U}^0(x))
\end{equation}
locally in time has the same smoothness as the initial data \cite{Dafermos}. However, over a finite time $t_*$
the derivatives of solution can go to infinity, which corresponds to formation of singularities of the solution.
 The exact class of smooth initial data (\ref{CD}), leading to the formation of singularities of a solution to the system (\ref{e1})-- (\ref{e_state}) is currently unknown.
There are only partial results that allow us to find sufficient conditions for the initial data for which during a finite time, a singularity of the solution arises \cite{R98}, or the solution remains smooth for $t> 0$ \cite{R2005}.

Therefore,  there arises the question of models that may be close to the original, but they are easier to study. Such a model can be obtained, in particular, under the assumption of a small pressure. Namely, formally setting $ p = 0 $ in the equation (\ref {e1}), we obtain  the equations of gas dynamics "without pressure", a non-strictly hyperbolic system in which the equation for velocity (transport equation) is split off. The "pressureless" gas dynamics system is the simplest model of a highly rarefied medium, in particular,  intergalactic dust \cite{Shandarin}. It is characterized by the fact that delta-like singularities arise in the density component, in contrast to ordinary shock waves. The passage to the limit under the condition of low pressure was strictly justified in \cite{Chen} for  irrotational coordinate system.  Similar reasoning can be carried out in the presence of the Coriolis force. Thus, obtaining a solution of the vectorial transport equation
\begin{equation}\label{main}
 \partial_t {\bf U}+ ({\bf U},\nabla) {\bf U} + lL {\bf U}= 0,
\end{equation}
we can find the velocity and entropy components from linear equations
(\ref{e2}) and (\ref{e3}). This approach was, in particular, used in
\cite{LiuTadmor}, where equation (\ref{main}) was investigated. Namely, the authors found conditions  on the initial data under which the singularity is formed and  concluded that the rotation of the coordinate system prevents the formation of singularities.
Indeed, an implicit solution to the Cauchy problem (\ref{main}),
\begin{equation}\label{main_cd}
{\bf U}(t,x)|_{t=0}=  {\bf U}_0(x)
\end{equation}
is easy to find (see below). However, the implicit form gives little information about the behavior of solution, in particular, about the nature of singularities. Unlike previous works, will get an explicit asymptotic representation of the solution until the appearance of singularity, which allows us to understand what happens to the solution near the critical time. Our method also allows to calculate this critical time accurately. The main tool is the method of stochastic perturbation along characteristics used in  \cite{AR}, \cite{AKR}.

\section{Implicit representation of solution of the transport equation}

Let us denote $U_1=u,$ $U_2=v$, $x_1=x,$ $x_2=y.$  The Cauchy problem for
(\ref{main}) takes the form
\begin{equation}\label{main1}\frac{\partial u}{\partial t}+u \frac{\partial u}{\partial x}+v \frac{\partial u}{\partial y}-l v=0,
\qquad \frac{\partial v}{\partial t}+u \frac{\partial v}{\partial
x}+v \frac{\partial v}{\partial y}+l u=0,\end{equation}
\begin{equation}\label{CD1}(u (0,x,y), v (0,x,y))= (u_{0}(x,y) , v_{0}(x,y)) \in C^1({\mathbb R}^2).\end{equation}
The equations of characteristics are the following:
\begin{equation}\label{char}
\dfrac{dx}{dt}=u,\quad \dfrac{dy}{dt}=v,\quad\dfrac{du}{dt}=
lv,\quad\dfrac{du}{dt}=lu.
\end{equation}
Thus, one can find a complete set of functionally independent first integrals
\begin{equation*} y-\frac{u}{l}=C_1,\quad x+\frac{v}{l}=C_2,\quad u
\sin{lt}+v \cos{lt}=C_3,\quad \quad u \cos{lt}-v \sin{lt}=C_4.
\end{equation*}
The implicit form of solution $(u,v)$ of (\ref{main1}) is
\begin{eqnarray*}
u=u_0\left(\frac{1}{l} (v(1-\cos lt)-u \sin lt)+x,\frac{1}{l}((\cos
lt-1)-v \sin lt)+y  \right), \\ v=v_0\left(\frac{1}{l} (v(1-\cos
lt)-u \sin lt)+x,\frac{1}{l}((\cos lt-1)-v \sin lt)+y  \right).
\end{eqnarray*}
This does not allow us to imagine the behavior of the solution.

\section{Criterion of singularity formation}

The system (\ref {main}) is hyperbolic, Thus, the loss of smoothness by a solution can be related to the unboundedness of either the solution itself or its derivatives \cite{Dafermos}. As follows from
(\ref{char}), $u$ and $v$ remain boundeded until the intersection of characteristics, that is, while the solution remains smooth.
Thus, we must find the necessary and sufficient conditions for the derivatives to go to infinity, and thereby find the class of initial data for which the solution remains smooth for all $ t> 0. $

Assuming $ C^1 $ -- smoothness of the solution, we can
differentiate (\ref {main}) with respect to spatial variables, obtaining along the characteristics the Riccati equation
 for the 2x2 matrix $ Q = (q_ {ij}) $ with components $ q_{11} = u'_x $, $ q_{12} = u'_y $, $ q_{21} = v'_x $, $ q_{22} = v'_y $:
\begin{equation} \label{Q0}
\dfrac{dQ}{dt}=-Q^{2}-lLQ.
\end{equation}
For $\det Q \ne 0$,  the change
  ${\mathcal Q}=Q^{-1}$ 
 reduces \eqref{Q0} to 
$\dfrac{
d{\mathcal Q}}{dt}=E+l{\mathcal Q} L$
 \cite{Egorov}. The latter linear equation can be  solved in a standard way. Solution (\ref{Q0}) has the form
\begin{equation}Q(t)=\frac{1}{l\det{\mathcal Q}(t)}
\begin{pmatrix}
K_2 l\cos lt-K_1l\sin lt  & 1-lK_3\cos lt +K_4l\sin lt \\
-1-lK_1\cos lt-lK_2\sin lt & lK_4\cos lt+lK_3\sin lt
\end{pmatrix}, \label{qsol}
\end{equation}
where
\begin{equation}\label{det}
\det {\mathcal Q}(t)= (1+l^2(K_2 K_4-K_1 K_3)+(K_2+K_4)l\sin lt
+(K_1 -K_3)l\cos lt)/l^2,
\end{equation}
 $K_1=-\frac{1}{l}-\frac{(v_0)_x}{D_0}$,
$K_2=\frac{(u_0)_x}{D_0}$, $K_3=\frac{1}{l}-\frac{(u_0)_y}{D_0}$,
$K_4=\frac{(v_0)_y}{D_0}$, $D_0={\rm det}Q\Big|_{t=0}.$ Solution of
 (\ref{Q0}) for ${\rm det} Q=0 $ can be obtained from
(\ref{qsol})  by passage to the limit as $D_0\to 0$.


\begin{thm}{} Assume  that initial data \eqref{main_cd} (\eqref{CD1}) belong to the class $C^1({\mathbb R}^2)$. If at $t=0$  for any point $(x,y)\in
{\mathbb R}^2$ condition
\begin{equation}\label{ba}
\Delta(x,y)\equiv
((u_0)_x+(v_0)_y)^2-4(D_0)-2l((v_0)_x-(u_0)_y)-l^2<0
\end{equation}
holds, then
 the solution to the Cauchy problem \eqref{main}, \eqref{main_cd} 
 remains $C^1$ -- smooth for all  $t>0$.
 \end{thm}

For the {\bf proof} it suffices to calculate that if  for any point $ (x, y) \in {\mathbb R}^2 $ condition \eqref{ba} is satisfied at $ t = 0 $, then the denominator in (\ref{qsol}) does not vanish, so the solution remains bounded for all $ t> 0 $. Otherwise, the derivatives of solution go to infinity in a finite time.


\begin{cor}{} For any initial data \eqref {main_cd} from  $ C^1 ({\mathbb R}^2) $, for which the first derivatives are uniformly bounded, we can choose $ l $ so large that the solution to the Cauchy problem \eqref{main}, \eqref{main_cd} remains $ C^1 $ - smooth for all
$t> 0$.
\end{cor}

Indeed, an analysis of the expression $ \Delta $ shows that  condition \eqref{ba}, which guarantees a time-smooth global solution, will always be satisfied for sufficiently large $ l $.


Condition (\ref{ba}) can be rewritten as
 \begin{equation*}\label{ba1}
({\rm div} {\bf U}^0)^2-4J({\bf U}^0)-2l\xi-l^2<0,
 \end{equation*}
where $J({\bf U}^0)={\rm det}\big(\frac{\partial U^0_i}{\partial
x_j}\big)$, $i,j=1,2,$ $\xi=(U^0_2)_{x_1}-(U^0_1)_{x_2}$.
A similar result in other terms was obtained in
\cite{LiuTadmor}.


\section{Representation of a solution based on the method of stochastic perturbation of characteristics}

The idea of the method is the following:  instead of the variables $ (x, y, u, v) $, we will consider their stochastic analogues $ (X, Y, U, V) $,
obeying the system of stochastic differential equations (SDE) 
\begin{equation}\label{sdu}
dX=U dt+ \sigma dW_1,\quad dY=V dt+ \sigma dW_2, \quad  dU=lV
dt,\quad dV=-lU dt,
\end{equation}
$$
X(0)=x^0,\quad Y(0)=y^0,\quad U(0)=u^0,\quad V(0)=v^0,\quad t\ge 0,
$$
where the random variables $ (X (t), Y (t), U (t), V (t)) $ belong to the phase space $ \Omega \times {\mathbb R}^2, $ $ \, \Omega \subset {\mathbb R}^2, $ \, $ W_t = (W_k)_t $, $ k = 1, 2 $ is the standard two-dimensional Brownian process, $ \sigma $ is a positive constant.

The system (\ref{sdu}) differs from (\ref{char}) by the added
stochastic disturbance along particle trajectories.

The Kolmogorov-Fokker-Planck equation describing the probability density $ P = P (t, x, y, u, v) $ in the space of coordinates and velocities has the form
\begin{eqnarray}\label{FPK}
\frac{\partial P(t,x, y,u,v)}{\partial t }=\\\left[-u \frac{\partial
}{\partial x }-v \frac{\partial }{\partial y }\, + l v\frac{\partial
}{\partial u }\,-\,  l u\frac{\partial }{\partial v }\, +
\frac{1}{2}\, \sigma^2 \,\left(\frac{\partial^2 }{\partial
x^2}\,+\,\frac{\partial^2 }{\partial y^2}\right)\right]
P(t,x,y,u,v),
\end{eqnarray}
with initial data
$$P(0,x,y,u,v)=P_0(x,y,u,v).$$


Denote by $ (\hat u (t, x, y), \hat v (t, x, y)) $ the mathematical expectation of the velocity $ (U (t), V (t)) $ at a fixed time $ t $ and  fixed value of $ (X (t), Y (t)) $. Then
\begin{equation}\label{Expect}
\hat u(t,x,y)=\frac{\int\limits_{{\mathbb
R}^2}\,u\,P(t,x,y,u,v)\,du\, dv}{\int\limits_{{\mathbb
R}^2}\,P(t,x,y,u,v)\,du\,dv},\quad \hat
v(t,x,y)=\frac{\int\limits_{{\mathbb R}^2}\,v\,P(t,x,y,u,v)\,du\,
dv}{\int\limits_{{\mathbb R}^2}\,P(t,x,y,u,v)\,du\,dv},
\end{equation}
$t\ge
0,\,(x,y)\in \Omega$.
We choose \begin{equation}\label{CD_FPK}
 P_0(x,y,u,v)=\delta
(u-u_0(x,y))\,\delta (v-v_0(x,y))\,f_0(x,y)\end{equation}  with arbitrary sufficiently smooth function $f_0(x)$ to obtain
$$
\hat u(0,x,y)=u_0(x,y), \quad \hat v(0,x,y)=v_0(x,y).
$$
Function $f_0(x,y)$ has a sense of probability density of the particle position in space at the initial moment of time, therefore, if we want to connect the system (\ref{main}) with the continuity equation (\ref{e2}), then we should choose $ f_0 (x, y) = \rho^0 (x, y) $.

Denote by $ \lambda_1, \lambda_2, \xi_1, \xi_2 $ 
variables dual to $ x, y, u, v $ respectively. Equation for
Fourier transform $ \tilde P (t, \lambda_1, \lambda_2, \xi_1, \xi_2) $ of the solution of the equation (\ref{FPK}) has the form
\begin{equation}\label{FT_FPK}
\frac{\partial \tilde{P}}{\partial t}= -\frac{1}{2} \tilde{P} \sigma
^2 ( \lambda_1^2+\lambda_2^2) +(\lambda _1-l\xi_2)\frac{\partial
\tilde{P}}{\partial \xi _1}+ (\lambda _2+l \xi_1)\frac{\partial
\tilde{P}}{\partial \xi_2},
\end{equation}
the Cauchy data (\ref{CD_FPK}) are transformed as
\begin{equation}\label{FT_CD}
\tilde{P}(0,\lambda_1, \lambda_2 ,\xi_1, \xi_2)=\int _{{\mathbb
R}^2}e^{-i(\lambda ,s)}e^{-i\left(\xi_1 u_0(s)+\xi_2
v_0(s)\right)}f_0(s)ds, \quad s=(s_1,s_2).
\end{equation}
The solution to (\ref{FT_FPK}), (\ref{FT_CD}) is
\begin{eqnarray*}
\tilde{P}(t,\lambda_1, \lambda_2 ,\xi_1, \xi_2)=\phantom{\int_{{\mathbb R}^2} f_0(s)
e^{-i (\lambda,s)}e^{-i(v_0(s)(a\cos(lt)+b\sin(lt)+\lambda _1)+
u_0(s)(a\sin(lt)-b \cos(lt)-\lambda _2))}}\\e^{-\frac{1}{2}
\sigma ^2\left| \lambda \right| ^2 t} \int_{{\mathbb R}^2} f_0(s)
e^{-i (\lambda,s)}\,e^{-i(v_0(s)(a\cos(lt)+b\sin(lt)+\lambda _1)+
u_0(s)(a\sin(lt)-b \cos(lt)-\lambda _2))}ds,
\end{eqnarray*}
where $a=\xi _2-\lambda _1,$ $b=-\xi_1-\lambda _2$. Performing the inverse Fourier transform, we obtain
\begin{equation*}
P(t,x,y,u,v)=\frac{1}{2 \pi \sigma ^2 t}\int_{{\mathbb R}^2}
\delta\left(u-q_1(t,s))\right)\delta \left(v-q_2(t,s)\right)
G(t,s,x,y,u,v)
 f_0(s)d s,
\end{equation*}
where $$q_1(t,s)=u_0(s) \cos lt-v_0(s) \sin lt,\qquad q_2(t,s)=v_0(s)
\cos lt+u_0(s) \sin lt,$$
$$G(t,s,x,y,u,v)=\exp
\left(-\frac{\left(u-u_0(s)+s_2-y\right)^2+\left(-(v-v_0(s))+s_1-x\right)^2}{2
\sigma ^2 t}\right).$$

 According to (\ref{Expect}) we finally get 
\begin{eqnarray}\label{Exspect_UV}
\hat{u}(t,x,y)=\frac{\int_{{\mathbb R}^2}  f_0(s) q_1(t,s)H(t,s,x,y)
ds}{\int_{{\mathbb R}^2} f_0(s)H(t,s,x,y) ds},\,
\hat{v}(t,x,y)=\frac{\int_{{\mathbb R}^2}  f_0(s) q_2(t,s)H(t,s,x,y)
ds}{\int_{{\mathbb R}^2} f_0(s)H(t,s,x,y) ds},
\end{eqnarray}
where \begin{eqnarray*}H(t,s,x,y)=&\\&\exp 
\left(-\frac{\left(q_1(t,s_1,s_2)-u_0(s_1,s_2)+s_2-y\right)^2+\left(-q_2(t,s_1,s_2)+v_0(s_1,s_2))+s_1-x\right)^2}{2
\sigma ^2 t}\right).
\end{eqnarray*}

Completely analogously   to Proposition 1 (\cite{AKR}) one can  show that  if the functions $ u_0, $ $ v_0 $ and $ f_0> 0 $ belong to the class $ C^1 (\mathbb {R}^2) \cap C_b (\mathbb {R}^2), $ then as long as the solution to the Cauchy problem (\ref{main1}), (\ref{CD1}) remains smooth, the functions $ (\hat {u} (t, x, y ), \hat {v} (t, x, y) $ tend to solve the problem (\ref {main1}), (\ref {CD1}) as $ \sigma \rightarrow 0 $ for every fixed $ (t, x , y) $. This conclusion also follows from the results of \cite{Freidlin1}.

The following theorem sums up our reasoning.

\begin{thm} The classical solution of the Cauchy problem \eqref{main}, \eqref{main_cd} 
can be obtained by passing to the limit as $ \sigma \to 0 $ from the integral asymptotic representation \eqref{Exspect_UV}.
 \end{thm}

\section{An example of explicit representation of solution}

As the initial data, we choose a vortex solution
\begin{equation}\label{uv_data}
u_1\,=\,B_0\,\mu\,x_2\,e^{-\frac{\mu}{2}\, \left(
{x_1}^{2}+{x_2}^{2} \right)},\qquad u_2\,=\,
-B_0\,\mu\,x_1\,e^{-\frac{\mu}{2}\, \left( {x_1}^{2}+{x_2}^{2}
\right)},
\end{equation}
with  $B_0=\rm{const}$ and $\mu=\rm{const}>0$.  As shown in \cite {RYH}, this velocity field is part of the stationary solution $(\rho,{\bf
U})$ of system \eqref{e1} -- \eqref{e3} for constant $S$
and the state equation $p=c_0 \rho^\gamma$,  $c_0=\rm{const}>0$. In this case, the density
$\rho=\pi^\frac{2\gamma}{\gamma-1}$, where
\begin{equation*}\label{pi}
\pi(x_1,x_2)\,=\,-\,\frac{1}{c_0}\,\left(\frac{1}{2}\,B_0^2\,\mu\,e^{-\mu\,(x_1^2+x_2^2)}\,
-\,l\,B_0\,e^{-\frac{\mu}{2}\,(x_1^2+x_2^2)}\right).
\end{equation*}
In the case of  system (\ref {main}), (\ref {e2}) this solution, of course, will no longer be stationary, but it can be chosen as the initial data in order to study the formation of singularities.

We will be interested in the following questions:

1. At what point in time does the formation of a singularity begin?

2. At what point in space does the singularity first appear?

3. What is the type of emerging singularity?

\begin{figure}[h]
\begin{minipage}{0.3\columnwidth}
\centerline{\includegraphics[width=1.2\columnwidth]{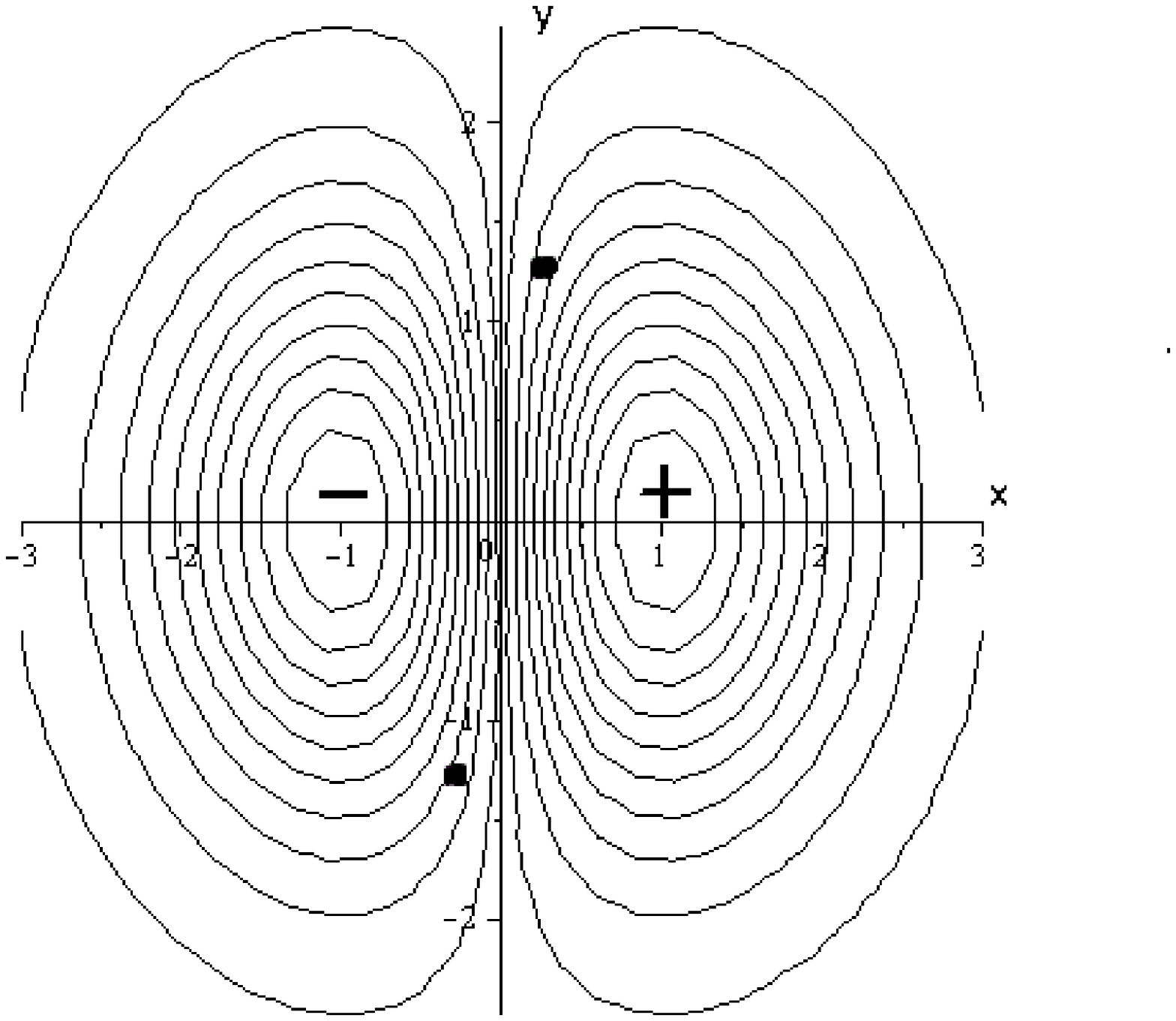}}
\end{minipage}
\begin{minipage}{0.3\columnwidth}
\centerline{\includegraphics[width=1.2\columnwidth]{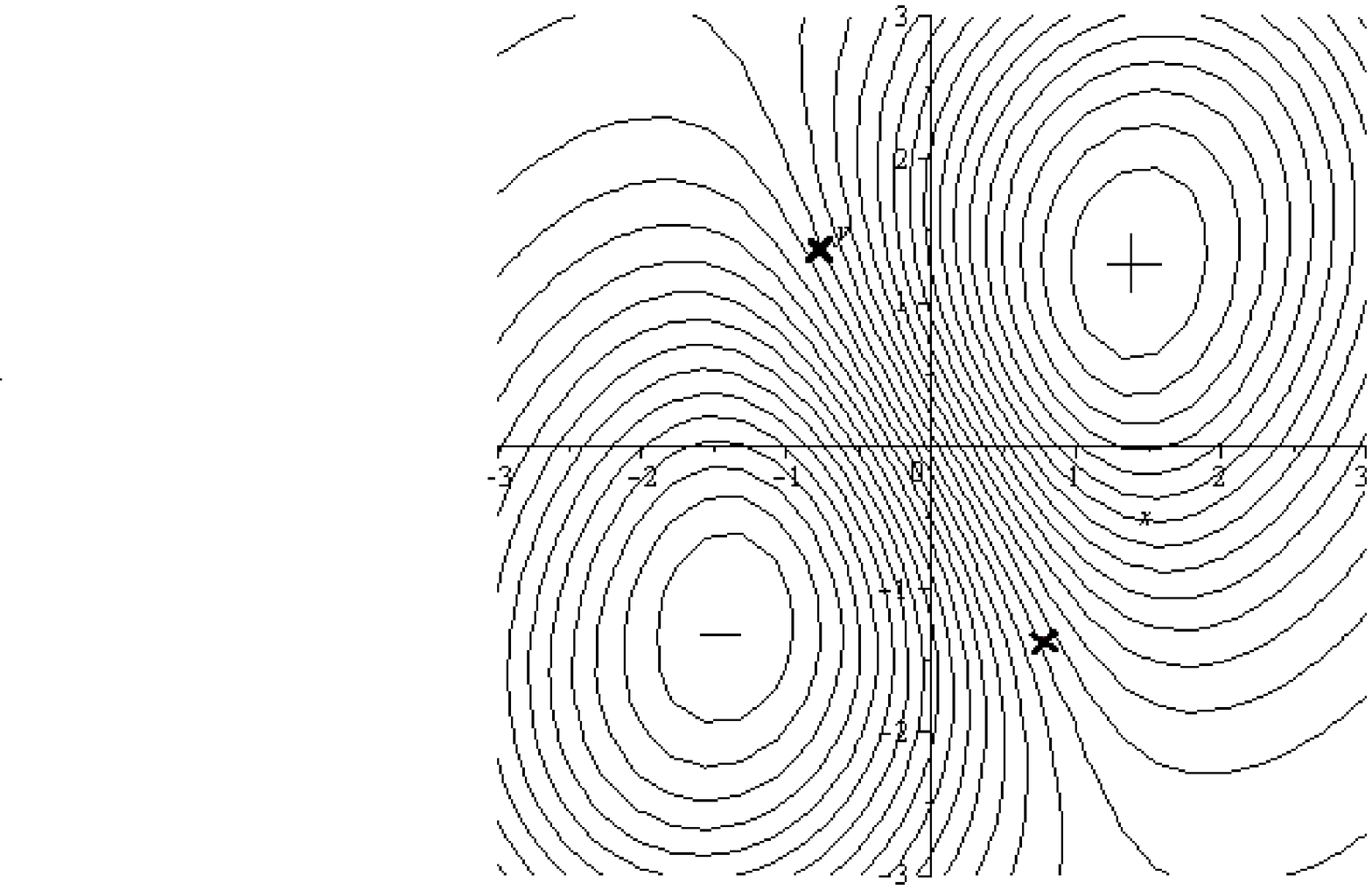}}
\end{minipage}
\begin{minipage}{0.3\columnwidth}
\centerline{\includegraphics[width=1.2\columnwidth]{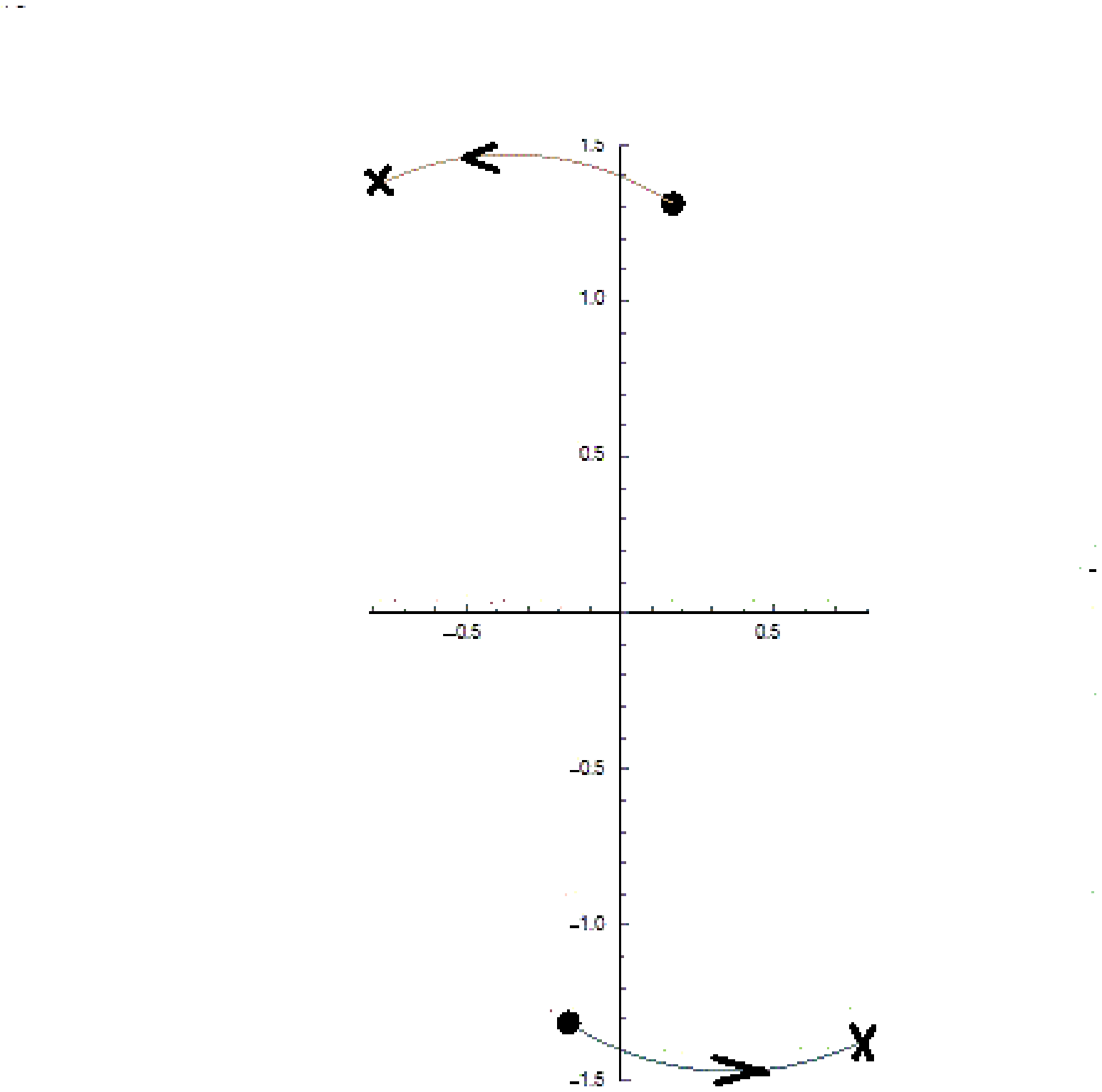}}
\end{minipage}
\caption{Level lines of  $ v $ at the initial moment of time (left); Level lines of $ \hat v $ 
at the moment of singularity formation (center);
Lagrangian trajectories leading to points of singularities formation
(right). }
\end{figure}

To get graphical illustrations, we choose the following  values of parameters: $ l = 1 $, $ \mu = 1 $, $ B_0 = -1.5 $. For simplicity, we will consider all quantities to be dimensionless. First, find the points of the plane $ {\mathbb R}^2 $ for which the determinant 
\eqref{det} vanishes first of all, that is, $$ (x_1, x_2) | \min \limits_{{\mathbb R}^2 } \Delta (x_1, x_2) \ge 0, $$ where $ \Delta $ is defined in \eqref {ba}. For the data \eqref {uv_data} this occurs at two points located symmetrically with respect to the origin, their position is marked in Fig.1 with dots.

We calculate the time $ t_* $ such that $ \det {\mathcal Q} (t) $ vanishes for selected points, $ t_* \approx 1.9 $. 
 Fig.1 (left) shows the level lines of the initial datum $ v_0 (x, y) $ and the position of the points that generate the singularities in the first place (for example, we chose one of the components of solution). Fig.1 (center) shows the level lines of the function $ \hat v(t, x, y)$ found by the formula \eqref{Exspect_UV} for $ \sigma = 0.1 $, $ t = t_* $ and the position of the points in which singularities arises. At these points, as follows from \eqref{qsol}, the derivatives of solution go to infinity. 
On the right, Fig.1 shows the Lagrangian trajectories found by  formulae \eqref{char} leading from the initial state to the points of singularities formation, marked with crosses.

Note that $ u (t, x, y) $ can be obtained from $ \hat u (t, x, y) $
when passing to the limit, therefore, the level lines give an idea of the behavior of the limit function only approximately. 
For very small $ \sigma $
drawing of $ \hat u (t, x, y) $  has difficulties associated with the numerical computation of improper integrals.

\section{Conclusion}

We study the classical solutions of the two-dimensional transport equation in the presence of the Coriolis force. The class of initial data for which the solution remains smooth for all $ t \ge 0 $ is found. An asymptotic representation of the solution is obtained until the moment of loss of smoothness. The process of singularities formation is considered for specific initial data.
Note that the representation of a solution based on the stochastic regularization method is related to the  vanishing viscosity method, but its scope is much wider. For example, even in the absence of the Coriolis force, using the Cole-Hopf transform, we can obtain an integral formula for representing the solution only in case of potential flow \cite {R_HJ}. The corresponding integral formula based on the stochastic regularization method can be obtained for arbitrary initial data in the presence of various external forces, including the Coriolis force.


\begin{thebibliography}{99}

\bibitem {Vallis}  Vallis, G.K.  Atmospheric and oceanic fluid dynamics. Fundamentals and large-scale
circulation. Cambridge University Press (2006).



\bibitem{Obukhov} Obukhov,  A.M. On the geostrophical wind. Izv.Acad.Nauk (Izvestiya of Academie of Science
of URSS), Ser. Geography and Geophysics, XIII, 281-306 (1949).

\bibitem{Dafermos} {Dafermos~C.M.}
Hyperbolic Conservation Laws in Continuum Physics. The 4th Edition,
Berlin-Heidelberg: Springer, 2016.

\bibitem{R98} Rozanova O.S. The formation of singularities of solutions with a compact support of Euler equations on a rotating plane. Differential equations,   34 (1998), N 8, 1114-1118.

\bibitem{R2005} Rozanova O.S. Classes of smooth solutions to multidimensional balance laws of gas dynamic
type on Riemannian manifolds. In: Trends in mathematical physics
research,  Nova Sci. Publ Hauppauge, NY, 155-204, 2004.

\bibitem{Shandarin}  Shandarin S.F.,  Zeldovich Y.B., Large-scale structure of the universe:
turbulence, intermittency, structures in a selfgravitating medium,
Rev. Modern Phys. 61 (1989) 185-220.


\bibitem{Chen} Chen G.-Q., Liu H. Concentration and cavitation in the vanishing
pressure limit of solutions to the Euler equations for nonisentropic
fluids. Physica D: Nonlinear Phenomena, 189 (2004),141-165.

\bibitem{LiuTadmor} H.Liu, E.Tadmor,  { Rotation prevents finite-time
breakdown} Physica D: Nonlinear Phenomena, { 188}(2004) 262--276.

\bibitem{AR}Korshunova A., Rozanova O. On effects of stochastic regularization for the pressureless gas dynamics
 Contemp. Appl. Math., 17 (2012), 486-493.

\bibitem{AKR}Albeveriio S., Korshunova A., Rozanova O.  A probabilistic model associated with the pressureless gas dynamics,
Bulletin des Sciences Mathematiques,  137 (2013), 902-922.





\bibitem{Egorov} Egorov, A.I. Riccati Equations.
Pensoft Publishers, 2007.

\bibitem {Freidlin1} M.Freidlin,\, {Markov processes and differential equations:
asymptotic problems}, Lectures in Mathematics, ETH Z\"urich. Basel:
Birkha\"user, 1996.

\bibitem {RYH} Rozanova, O.S.,   Yu, J.-L.,  Hu, C.-K.: On the position of vortex in
two-dimensional model of atmosphere. Nonlinear Analysis: Real World
Applications, 13(2012), 1941--1954.

\bibitem{R_HJ} Rozanova O.S. On the connection of the Hamilton-Jacobi equation and some systems of quasilinear equations.
Proceedings of the Institute of Mathematics, Ural Branch of the Russian Academy of Sciences, 21 (2015),  N2, 206-219.
\end{thebibliography}
\end{document}